\documentclass[11pt,a4paper]{article}

\usepackage[T1]{fontenc}
\usepackage[utf8]{inputenc}
\usepackage{lmodern}
\usepackage{microtype}
\usepackage{geometry}
\usepackage{amsmath,amssymb,amsthm,mathtools}
\usepackage{booktabs,longtable,array}
\usepackage{enumitem}
\usepackage[hidelinks]{hyperref}

\newtheorem{theorem}{Theorem}[section]
\newtheorem{proposition}[theorem]{Proposition}
\newtheorem{lemma}[theorem]{Lemma}
\newtheorem{corollary}[theorem]{Corollary}
\newtheorem{remark}[theorem]{Remark}
\newtheorem{definition}[theorem]{Definition}

\newcommand{\LimQ}{\mathcal L_{L,Q}^{\mathrm{bip}}}
\newcommand{\rh}{\rho}
\newcommand{\Ga}{\Gamma_L}
\newcommand{\bb}{\mathbf b}
\newcommand{\DD}{\mathcal D}

\title{Intervals of limit points for $B_{\alpha}$-spectral radii}
\author{%
  Francesco Belardo\textsuperscript{1}\quad
       Callum Huntington\textsuperscript{1}\\[0.35em]
  Elismar R. Oliveira\textsuperscript{2}\quad
       Vilmar Trevisan\textsuperscript{2}\\[0.75em]
  {\small\textsuperscript{1}Department of Mathematics and Applications
    ``R. Caccioppoli''}\\
  {\small University of Naples Federico II, Naples, Italy}\\
  {\small\texttt{fbelardo@unina.it},\quad
    \texttt{callum.huntington@unina.it}}\\[0.55em]
  {\small\textsuperscript{2}Instituto de Matem\'atica e Estat\'istica}\\
  {\small Universidade Federal do Rio Grande do Sul, Porto Alegre, Brazil}\\
  {\small\texttt{elismar.oliveira@ufrgs.br},\quad
    \texttt{trevisan@mat.ufrgs.br}}%
}
\date{}
\begin{document}
\maketitle
\begin{abstract}
For a finite simple graph $G$, let
\[
 B_{\alpha}(G)=\alpha A(G)+(1-\alpha)L(G)
       =(1-\alpha)D(G)+(2\alpha-1)A(G),
       \qquad 0\leq \alpha\leq 1.
\]
We study the limit points of the spectral radii of these matrices.  Put
$a_{\alpha}=1-\alpha$, $b_{\alpha}=|2\alpha-1|$, and
$w_{\alpha}=a_{\alpha}+b_{\alpha}$.  We prove that, for every
$\alpha\neq 1/2$ and every integer $d\geq 3$, the interval
\[
 \left[a_{\alpha}d+2b_{\alpha}\sqrt{d-1},\,w_{\alpha}d\right]
\]
consists entirely of $B_{\alpha}$-limit points.  All these points are
realized by connected bipartite graphs of maximum degree at most $d$.
We extend the induced-subgraph interpolation of Alon and Wei from $A$ to
$aD+bA$, and obtain a unified proof using bipartite Ramanujan graphs.
Consecutive intervals overlap
for all sufficiently large $d$; hence the set of $B_{\alpha}$-limit points
contains a half-line for every $\alpha\neq 1/2$.  At the exceptional value
$\alpha=1/2$, one has $B_{1/2}=D/2$, and the set of spectral radii is discrete.
Combining uniform edge subdivisions with generalized Bethe decorations,
we improve the Laplacian and signless Laplacian half-line to $[29/5,\infty)$.
The finite covering step is verified by an exact rational certificate;
all 101 intervals and the verification procedure are supplied.

We also treat two complementary tree constructions.  For the starlike trees
$T_{1,m,m}$, the limiting $B_{\alpha}$-spectral radius is the unique solution
above the path threshold of an explicit unsquared equation.  More generally,
letting selected arms of a starlike tree tend to infinity produces ordered limit points with tree realizations.  Finally, we give the exact affine relations among $B_{\alpha}$,
the $A_{\beta}$-matrix, and the deformed Laplacian.  These identities identify
precisely which tree and caterpillar results are reparametrizations of earlier
work.  The interval theorem also yields half-lines of $A_\beta$-limit
points for every $0\le\beta<1$, and of deformed-Laplacian limit points for
every fixed nonzero real deformation parameter.
\end{abstract}

\section{Introduction}

Let $A(G)$, $D(G)$, and $L(G)=D(G)-A(G)$ denote, respectively, the
adjacency, degree, and Laplacian matrices of a finite simple graph $G$.
Samanta, Deepshikha, and Das introduced the one-parameter family
\begin{equation}\label{eq:B-definition}
 B_{\alpha}(G)=\alpha A(G)+(1-\alpha)L(G)
 =(1-\alpha)D(G)+(2\alpha-1)A(G),
 \qquad 0\leq\alpha\leq1.
\end{equation}
It contains, up to a positive scalar, four standard graph matrices:
\[
 B_0=L,\qquad B_{1/2}=\tfrac12D,\qquad
 B_{2/3}=\tfrac13Q,\qquad B_1=A.
\]
Although $B_{\alpha}$ need not be entrywise nonnegative, its spectral radius
is its largest eigenvalue \cite{SamantaDeepshikhaDas}.  We write this value as
$\rho_{B_\alpha}(G)$.  Indeed, for $\alpha\le1/2$,
$B_\alpha=(1-2\alpha)L+\alpha D$ is positive semidefinite; for
$\alpha\ge1/2$, it is entrywise nonnegative and Perron--Frobenius applies.

The Hoffman program asks for the accumulation points of spectral radii of
graph matrices.  The adjacency case was initiated by Hoffman
\cite{Hoffman1972}; Shearer subsequently proved that every real number at
least $\sqrt{2+\sqrt5}$ is an adjacency limit point \cite{Shearer1989}.
Several variants have since been considered, including the Laplacian and
signless Laplacian matrices \cite{BelardoOliveiraTrevisan,Guo2008}, signed
adjacency matrices, and the $A_{\beta}$-matrix \cite{Nikiforov}
\[
 A_{\beta}(G)=\beta D(G)+(1-\beta)A(G)
 \qquad (0\leq\beta\leq1).
\]
Oliveira and Trevisan recently adapted Shearer's caterpillars to
$A_{\beta}$ for $\beta<1/2$ \cite{OliveiraTrevisanAalpha}.
The deformed Laplacian
\[
 M_G(s)=I-sA(G)+s^2(D(G)-I)
\]
has also been studied from the same limit-point perspective
\cite{Diaz2026spectral,Oliveira2027deformed}.
On bipartite graphs, switching and affine changes of parameters relate
these matrix families; the precise domains and the distinction between
largest eigenvalues and spectral radii are recorded in
Section~\ref{sec:relations}.

We use the strict, topological notion of a limit point.

\begin{definition}\label{def:limit-point}
For fixed $\alpha$, a real number $\lambda$ is a
$B_{\alpha}$-\emph{limit point} if there is a sequence of finite simple graphs
$(G_k)$ such that the numbers $\rho_{B_\alpha}(G_k)$ are pairwise distinct and
\[
 \lim_{k\to\infty}\rho_{B_\alpha}(G_k)=\lambda.
\]
The set of all such numbers is denoted by $\mathcal L_{B}(\alpha)$.
Since the spectral radius of a disconnected graph is attained by one of its
connected components, the graphs in the definition may be taken connected.
We write $\mathcal L_B^{\mathrm{bip}}(\alpha)$ when only connected
bipartite graphs are allowed.
\end{definition}

The main result identifies explicit intervals contained in
$\mathcal L_B(\alpha)$.  Throughout the paper we use
\begin{equation}\label{eq:parameters}
 a_{\alpha}:=1-\alpha,\qquad
 c_{\alpha}:=2\alpha-1,\qquad
 b_{\alpha}:=|c_{\alpha}|,\qquad
 w_{\alpha}:=a_{\alpha}+b_{\alpha}.
\end{equation}

\begin{theorem}[Main interval theorem]\label{thm:main-interval}
Let $\alpha\in[0,1]\setminus\{1/2\}$ and let $d\geq3$ be an integer.  Then
\begin{equation}\label{eq:Id}
 I_d(\alpha):=
 \left[a_{\alpha}d+2b_{\alpha}\sqrt{d-1},\,
       w_{\alpha}d\right]
 \subseteq \mathcal L_B(\alpha).
\end{equation}
More precisely, every point of $I_d(\alpha)$ is obtained as a limit of
pairwise distinct $B_\alpha$-spectral radii of connected bipartite graphs
with maximum degree at most $d$.
\end{theorem}

The proof is a weighted extension of the induced-subgraph interpolation in
\cite[Lemma~3.1]{AlonWei}; that interpolation strategy is already present
there for adjacency matrices.  We prove a sharper one-vertex averaging
bound for $aD+bA$, including the change in the degree diagonal.  Starting
from the bipartite Ramanujan graphs of Marcus, Spielman, and Srivastava
\cite{MSS}, we obtain the same interval for every $\alpha\ne1/2$ by
bipartite switching.  A short remark also records the regular-graph consequence of
\cite[Theorem~1.2]{AlonWei} for $\alpha<1/2$.

The intervals in \eqref{eq:Id} eventually overlap.   More precisely, for
$\alpha\neq1/2$ define
\begin{equation}\label{eq:dalpha-gamma}
 d_{\alpha}:=
 \left\lceil\left(1+\sqrt{\frac{w_{\alpha}}{b_{\alpha}}}\right)^2\right\rceil,
 \qquad
 \Gamma(\alpha):=
 a_{\alpha}d_{\alpha}
 +2b_{\alpha}\sqrt{d_{\alpha}-1}.
\end{equation}

\begin{corollary}[A half-line of limit points]\label{cor:halfline}
For every $\alpha\in[0,1]\setminus\{1/2\}$,
\[
 [\Gamma(\alpha),\infty)\subseteq\mathcal L_B^{\mathrm{bip}}(\alpha)\subseteq\mathcal L_B(\alpha).
\]
For $\alpha=1/2$, the set of possible spectral radii is
$\{k/2:k\in\mathbb N_0\}$ and consequently
$\mathcal L_B(1/2)=\varnothing$.
\end{corollary}

The weighted theorem alone gives the common Laplacian and signless
Laplacian threshold $\Ga=6+2\sqrt5$. Our second result lowers it by
combining subdivisions and uniform finite rooted-tree decorations.

\begin{theorem}[Laplacian and signless Laplacian half-line]\label{cor:LQ}
Every real number $\lambda\ge29/5$ is a limit of pairwise distinct
Laplacian spectral radii of connected bipartite graphs. The same graphs
give the same signless Laplacian spectral radii.
\end{theorem}

The proof is in Section~\ref{sec:bethe}. Its analytic part transfers
whole intervals through two graph operations; a chain of 101 such
intervals then covers $[29/5,\Ga]$. Appendix~\ref{app:chain} lists the
chain and explains its exact rational verification. The approximating
graphs may have cycles, so this is not a result restricted to trees.

For context, Guo's classical Laplacian sequence
$\beta_0=4$, $\beta_n=\lim_{m\to\infty}\rho_L(T_{1,n,m})$ for $n\ge1$
increases to $\mathfrak h_L=4.382975767\ldots$ \cite{Guo2008}.
Subdivision gives further nontrivial intervals tending to $9/2$.
These inclusions do not classify the full limit-point set, nor do they
establish full coverage of $[9/2,29/5)$.

Section~\ref{sec:unified} proves the weighted interval theorem and its
overlap corollary before the applications. Section~\ref{sec:bethe}
establishes the improved Laplacian half-line. Section~\ref{sec:trees}
collects the complementary tree constructions, while
Section~\ref{sec:relations} gives the exact relations with $A_\beta$
and the deformed Laplacian, including the results inherited from earlier
work.

\section{A unified interval theorem on bipartite graphs}\label{sec:unified}

For $a\ge0$ and $b>0$, write
\[
H_{a,b}(G)=aD(G)+bA(G),\qquad
\rho_{a,b}(G)=\lambda_1(H_{a,b}(G)).
\]
These matrices are nonnegative, so $\rho_{a,b}$ is also their spectral
radius.  The proof extends the induced-subgraph interpolation of
\cite[Lemma~3.1]{AlonWei} from $A$ to $aD+bA$.  The degree diagonal must
be recomputed after each deletion; Lemma~\ref{lem:deletion} includes this
change explicitly through an edge-counting identity.

Equation \eqref{eq:B-definition} can be written as
\begin{equation}\label{eq:B-ac}
 B_{\alpha}(G)=a_{\alpha}D(G)+c_{\alpha}A(G).
\end{equation}
The following observations will be used repeatedly.

\begin{lemma}[Bipartite switching]\label{lem:bipartite-switching}
Let $G$ be bipartite.  Then $B_{\alpha}(G)$ is similar to
\[
 H_{\alpha}(G):=a_{\alpha}D(G)+b_{\alpha}A(G).
\]
In particular,
\begin{equation}\label{eq:bipartite-A-beta}
 \rho_{B_\alpha}(G)
 =w_{\alpha}\rho_{A_{\gamma_\alpha}}(G),
 \qquad
 \gamma_{\alpha}:=\frac{a_{\alpha}}{w_{\alpha}}.
\end{equation}
\end{lemma}

\begin{proof}
If $c_{\alpha}\geq0$, the assertion is immediate.  If
$c_{\alpha}<0$, order the vertices by a bipartition $X\cup Y$ and put
\[
A=\begin{pmatrix}0&C\\C^T&0\end{pmatrix},\qquad
S=\begin{pmatrix}I_{|X|}&0\\0&-I_{|Y|}\end{pmatrix}.
\]
Then $S^{-1}=S$, $SDS=D$, and $SAS=-A$.  Hence this is an orthogonal
similarity, not merely a congruence, and
\[
 SB_{\alpha}(G)S=a_{\alpha}D(G)-c_{\alpha}A(G)
 =a_{\alpha}D(G)+b_{\alpha}A(G).
\]
Finally,
$a_{\alpha}D+b_{\alpha}A=w_{\alpha}A_{\gamma_\alpha}$.
\end{proof}

For $\alpha\geq1/2$ no bipartiteness assumption is needed:
\begin{equation}\label{eq:global-A-beta}
 B_{\alpha}(G)=\alpha A_{\beta}(G),
 \qquad
 \beta=\frac{1-\alpha}{\alpha}.
\end{equation}
Indeed, $\alpha\beta=1-\alpha$ and
$\alpha(1-\beta)=2\alpha-1$.  The range $\alpha>2/3$ corresponds exactly to $\beta<1/2$, the
range of the caterpillar half-line theorem of Oliveira and Trevisan \cite{OliveiraTrevisanAalpha}.

\subsection{Ramanujan graphs and two elementary estimates}
A connected $d$-regular bipartite graph is Ramanujan if all adjacency
eigenvalues other than the two trivial values $d$ and $-d$ have absolute
value at most $2\sqrt{d-1}$.  In particular its second largest adjacency
eigenvalue satisfies $\lambda_2\le2\sqrt{d-1}$; the negative eigenvalue
$-d$ creates no difficulty in an upper Rayleigh estimate.

\begin{theorem}[Marcus--Spielman--Srivastava, Theorem 5.5 of \cite{MSS}]
\label{thm:MSS}
For every integer $d\ge3$ there is a sequence of finite simple connected
$d$-regular bipartite Ramanujan graphs whose orders tend to infinity.
\end{theorem}

The construction starts from $K_{d,d}$ and repeatedly takes suitable
$2$-lifts.  Simplicity and bipartiteness are preserved.  The bound on the
new eigenvalues is strictly below $d$, so no additional eigenvalue $d$
appears and the lifts remain connected.
One could instead use near-Ramanujan graphs with
$\lambda_2\le2\sqrt{d-1}+o(1)$, as provided by Friedman \cite{Friedman}.
Control of merely a proportion tending to one of the nontrivial
eigenvalues would not suffice: the estimate requires control of $\lambda_2$.

\begin{lemma}[Small induced subgraphs]\label{lem:small-induced}
Let $G$ be a $d$-regular graph on $N$ vertices, with
$\lambda_2(G)\le r<d$.  For an induced subgraph $H=G[S]$ on $m$ vertices,
\begin{align}
\lambda_1(A(H))&\le r+(d-r)\frac mN,\label{eq:induced-adj-bound}\\
\rho_{a,b}(H)&\le ad+b\left(r+(d-r)\frac mN\right).
\label{eq:induced-M-bound}
\end{align}
\end{lemma}
\begin{proof}
Extend a unit vector on $S$ by zeros and write $x=tu+z$, where
$u=\mathbf1/\sqrt N$ and $z\perp u$.  Cauchy--Schwarz gives
$t^2\le m/N$, and therefore
\[
x^TA(G)x\le dt^2+r(1-t^2)\le r+(d-r)m/N.
\]
Maximizing proves the first inequality.  Since the degrees in $H$ are at
most $d$, $D(H)\le dI$, which proves the second.
\end{proof}

\begin{lemma}[Averaged vertex deletion]\label{lem:deletion}
Let $H$ have $k\ge2$ vertices and maximum degree at most $d$.  Some
vertex $v$ satisfies
\begin{equation}\label{eq:deletion-step}
0\le\rho_{a,b}(H)-\rho_{a,b}(H-v)
\le\frac{\rho_{a,b}(H)}{k-1}
\le\frac{(a+b)d}{k-1}.
\end{equation}
\end{lemma}
\begin{proof}
Put $\rho=\rho_{a,b}(H)$, choose a unit eigenvector $x$ for $\rho$, and
let $x^{(v)}$ be its restriction to $H-v$.  Then
\begin{equation}\label{eq:deletion-identities}
\sum_v\|x^{(v)}\|^2=k-1,\qquad
\sum_v(x^{(v)})^TH_{a,b}(H-v)x^{(v)}=(k-2)\rho.
\end{equation}
For the second identity, write
\[
x^TH_{a,b}(H)x
=\sum_{\{i,j\}\in E(H)}\bigl(a(x_i^2+x_j^2)+2bx_ix_j\bigr).
\]
Each edge contribution survives exactly $k-2$ vertex deletions.  This
accounts for both the adjacency terms and the decreases in all affected
degrees.  With $R=\max_v\rho_{a,b}(H-v)$, Rayleigh's inequality gives
$(k-2)\rho\le(k-1)R$, including zero restriction vectors without division.
Choose $v$ attaining $R$.  Monotonicity follows because $H_{a,b}(H-v)$
is entrywise no larger than the corresponding principal submatrix of
$H_{a,b}(H)$.  Finally, $\rho\le(a+b)d$ by the maximum row sum.
\end{proof}

\begin{lemma}[Strict accumulation from dense approximation]
\label{lem:strict-accumulation}
Let $E\subset\mathbb R$, and let $L<U$.  If $(L,U)\subset\overline E$,
then every point of $[L,U]$ is an accumulation point of $E$ in the sense
of a limit of pairwise distinct elements.
\end{lemma}
\begin{proof}
For a prescribed $x\in[L,U]$, choose points $t_j\in(L,U)$ tending
monotonically to $x$, and pairwise disjoint open intervals $V_j$ about
$t_j$, contained in $(L,U)$, with diameters tending to zero.  Since
$t_j\in\overline E$, choose $e_j\in E\cap V_j$.  Then the $e_j$ are
pairwise distinct and tend to $x$.
\end{proof}

\subsection{Interpolation and the main theorem}
\begin{theorem}[Positive-coefficient interpolation]\label{thm:positive-interpolation}
For $a\ge0$, $b>0$, and integer $d\ge3$, every point of
\begin{equation}\label{eq:Jab}
J_d(a,b)=[ad+2b\sqrt{d-1},(a+b)d]
\end{equation}
is a limit of pairwise distinct values $\rho_{a,b}(H)$ on connected
bipartite graphs $H$ of maximum degree at most $d$.
\end{theorem}
\begin{proof}
Put $r=2\sqrt{d-1}<d$, $L=ad+br$, and $U=(a+b)d$.
Take a graph $G$ of order $N$ from Theorem~\ref{thm:MSS} and set
$m=\lfloor\sqrt N\rfloor$.  Repeatedly apply Lemma~\ref{lem:deletion}
to obtain induced subgraphs
\[
G=H_N\supset H_{N-1}\supset\cdots\supset H_m,
\qquad |V(H_k)|=k.
\]
They remain bipartite and have maximum degree at most $d$.
Writing $r_k=\rho_{a,b}(H_k)$, the two estimates give
\[
r_N=U,\qquad
r_m\le L+b(d-r)m/N,
\qquad
0\le r_k-r_{k-1}\le\frac{(a+b)d}{m}\quad(k>m).
\]
Fix $t\in(L,U)$.  For sufficiently large $N$, $r_m<t<r_N$.
Hence there is an index $k$ with $r_{k-1}\le t\le r_k$, and
the crossing error is at most $(a+b)d/m=O(N^{-1/2})$.
Thus every interior target is in the closure of the set of these graph
indices.  Any disconnected crossing graph can be replaced by a connected
component attaining its index, preserving bipartiteness and the degree
bound.  Lemma~\ref{lem:strict-accumulation} supplies pairwise distinct
values, including at both endpoints.
\end{proof}

\begin{proof}[Proof of Theorem~\ref{thm:main-interval}]
Apply Theorem~\ref{thm:positive-interpolation} with
$a=a_\alpha$ and $b=b_\alpha>0$.  Every approximating graph is bipartite,
so Lemma~\ref{lem:bipartite-switching} identifies its $\rho_{a,b}$ with
$\rho_{B_\alpha}$, for either sign of $c_\alpha$.  The interval
\eqref{eq:Jab} is exactly $I_d(\alpha)$.
\end{proof}

\begin{remark}[A quantitative approximation]\label{rem:rate}
For every $t\in J_d(a,b)$, the same chain contains an index at distance
at most
\[
\max\left\{b(d-2\sqrt{d-1})\frac mN,\frac{(a+b)d}{m}\right\}
=O(N^{-1/2})
\]
from $t$: use a crossing when one exists, and $H_m$ otherwise.
At the upper endpoint $G$ itself has index $U$; strict accumulation at
that endpoint is established separately by the preceding lemma.  The
construction gives an existence bound; it does not bound the running time
of constructing the Ramanujan starting graphs.
\end{remark}

\subsection{The overlap argument}\label{sec:overlap}

We now prove Corollary~\ref{cor:halfline}.  Assume
$\alpha\neq1/2$, abbreviate $a=a_{\alpha}$, $b=b_{\alpha}$, and
$w=a+b$, and compare two consecutive intervals.  Their overlap condition is
\begin{align*}
 a(d+1)+2b\sqrt d&\leq wd\\
 &\Longleftrightarrow
 b(d-2\sqrt d)\geq a.                         \tag{\(*\)}\label{eq:overlap}
\end{align*}
The function $d\mapsto d-2\sqrt d$ is increasing for $d>1$.  Solving
\eqref{eq:overlap} for real $d$ gives
\[
 \sqrt d\geq1+\sqrt{1+a/b}
 =1+\sqrt{w/b}.
\]
Consequently, the integer $d_{\alpha}$ in
\eqref{eq:dalpha-gamma} is the first degree from which all consecutive
intervals overlap.  The union of $I_d(\alpha)$ for $d\geq d_{\alpha}$ is
therefore
\[
 [a d_{\alpha}+2b\sqrt{d_{\alpha}-1},\infty)
 =[\Gamma(\alpha),\infty).
\]
This proves the first assertion of Corollary~\ref{cor:halfline}, with bipartite realizations.

If $\alpha=1/2$, then
\[
 \rho_{B_{1/2}}(G)=\tfrac12\Delta(G).
\]
The set of all possible values is the discrete set
$\{k/2:k\in\mathbb N_0\}$, which has no finite accumulation point.  This
proves the exceptional assertion.  For $a=b=1$, the intervals
$[d+2\sqrt{d-1},2d]$ overlap for $d\ge6$, and hence
\begin{equation}\label{eq:oldray}
[\Ga,\infty)\subseteq\LimQ,\qquad \Ga=6+2\sqrt5,
\end{equation}
where $\LimQ$ denotes the common strict limit-point set of $L$ and $Q$
on connected bipartite graphs. This is the starting half-line for the
improvement in Theorem~\ref{cor:LQ}.

\begin{remark}[The first uncovered interval]\label{rem:first-gap}
The left endpoint of $I_4(\alpha)$ exceeds the right endpoint of
$I_3(\alpha)$ by $a_\alpha+b_\alpha(2\sqrt3-3)>0$.
Thus $I_3$ and $I_4$ never overlap.  The pair $I_4,I_5$ overlaps only
when $a_\alpha=0$, that is, at $\alpha=1$, where they meet at $4$.
For every other nonexceptional parameter $d_\alpha\ge5$.
These are gaps in the coverage by the displayed intervals, not assertions
that the full limit-point set has gaps there; for adjacency Shearer's
theorem already fills this part of the line.
\end{remark}

Notice that $\Gamma(\alpha)\to\infty$ as
$\alpha\to1/2$ from either side.  Thus the divergence of the sufficient
threshold is consistent with the genuinely discrete behavior at
$\alpha=1/2$.

\begin{remark}[The regular-graph route]\label{sec:regular-negative}
For $\alpha<1/2$ and a $d$-regular graph,
$\rho_{B_\alpha}(G)=a_\alpha d-b_\alpha\lambda_n(A(G))$.
Alon and Wei \cite[Theorem~1.2]{AlonWei} show that the least adjacency
eigenvalues of growing $d$-regular graphs with growing odd girth have
limit set $[-d,-2\sqrt{d-1}]$. This gives $I_d(\alpha)$ by the displayed
affine map; Lemma~\ref{lem:strict-accumulation} supplies pairwise distinct
values, including at the endpoints. A component attaining the index
remains $d$-regular. The weighted proof above additionally supplies
bipartite realizations and treats both signs of $2\alpha-1$ at once.
\end{remark}

\section{Subdivision and generalized Bethe decorations}\label{sec:bethe}
\subsection{Subdivision intervals}
Let $S_\ell(G)$ be obtained by replacing every edge of $G$ by a path of
exactly $\ell$ edges, where $\ell\ge1$. Define
\[
 p_{-1}(x)=0,\quad p_0(x)=1,\quad
 p_j(x)=(x-2)p_{j-1}(x)-p_{j-2}(x)\quad(j\ge1),
\]
and, for $x>4$,
\begin{equation}\label{eq:ab}
 a_\ell(x)=1+\frac{p_{\ell-2}(x)}{p_{\ell-1}(x)},\qquad
 b_\ell(x)=\frac1{p_{\ell-1}(x)}.
\end{equation}
For $\ell=1$ these are both $1$.

\begin{lemma}[Schur complement]\label{lem:schur}
For $x>4$,
\begin{equation}\label{eq:schur}
 \rh(Q(S_\ell(G)))<x
 \quad\Longleftrightarrow\quad
 \rh(a_\ell(x)D(G)+b_\ell(x)A(G))<x.
\end{equation}
The same equivalence holds with equality in both places.
\end{lemma}
\begin{proof}
The claim is immediate for $\ell=1$. For $\ell\ge2$, the block of
$xI-Q(S_\ell(G))$ on the internal vertices of each subdivided edge is
\[
 C=(x-2)I-A(P_{\ell-1}),
\]
which is positive definite for $x>4$. Its inverse has corner entries
$(C^{-1})_{11}=p_{\ell-2}/p_{\ell-1}$ and
$(C^{-1})_{1,\ell-1}=1/p_{\ell-1}$. Eliminating all these blocks gives
the Schur complement
\[
               xI-a_\ell(x)D(G)-b_\ell(x)A(G).
\]
Positive definiteness proves the strict inequality equivalence. Positive
semidefiniteness and singularity prove the equality equivalence.
\end{proof}

\begin{theorem}[Subdivision intervals]\label{thm:sub}
For $d\ge3$ and $\ell\ge1$, define $J_{d,\ell}=[L_{d,\ell},R_{d,\ell}]$,
where the endpoints are the unique solutions above $4$ of
\begin{align}
 (x-d)p_{\ell-1}(x)-dp_{\ell-2}(x)&=2\sqrt{d-1}
                  &&\text{for }L_{d,\ell},\label{eq:left}\\
 (x-d)p_{\ell-1}(x)-dp_{\ell-2}(x)&=d
                  &&\text{for }R_{d,\ell}.\label{eq:right}
\end{align}
Then $J_{d,\ell}$ is a nontrivial interval contained in $\LimQ$.
\end{theorem}
\begin{proof}
The equations are equivalent to $x=da_\ell(x)+cb_\ell(x)$ with
$c=2\sqrt{d-1}$ or $c=d$. For $\ell\ge2$, both $a_\ell-1$ and $b_\ell$
are nonnegative and decreasing: they are entries of the positive resolvent
$C^{-1}$, whose derivative is $-C^{-2}$. For $\ell=1$ they are constant.
Thus $x-da_\ell(x)-cb_\ell(x)$ increases at least as fast as $x$.
At $x=4$ it equals $4-2d+(d-c)/\ell<0$, whereas it tends to infinity.
The roots exist uniquely, and $2\sqrt{d-1}<d$ makes $L_{d,\ell}<R_{d,\ell}$.

Fix $x\in J_{d,\ell}$ and put $a=a_\ell(x)$, $b=b_\ell(x)$. The endpoint
characterization implies
\[
                ad+2b\sqrt{d-1}\le x\le(a+b)d.
\]
Theorem~\ref{thm:positive-interpolation} gives connected bipartite $G_n$ with pairwise
distinct values $\mu_n=\rh(aD(G_n)+bA(G_n))\to x$. Discard a possible
single value equal to $x$. Set
\[
 f_n(y)=y-\rh(a_\ell(y)D(G_n)+b_\ell(y)A(G_n))\qquad(y>4).
\]
Entrywise monotonicity gives $f_n(v)-f_n(u)\ge v-u$ for $4<u<v$.
Since $f_n(x)=x-\mu_n\to0$, for each $0<\varepsilon<x-4$ and all large
$n$ we have $f_n(x-\varepsilon)<0<f_n(x+\varepsilon)$. By
Lemma~\ref{lem:schur}, the radii of $S_\ell(G_n)$ tend to $x$. None equals
$x$, again by that lemma; a subsequence has pairwise distinct radii.
Subdivision preserves connectedness and bipartiteness.
\end{proof}

In particular, $J_{d,1}=I_d(0)$ and
\begin{equation}\label{eq:ell2}
 J_{d,2}=\left[
 \frac{d+2+\sqrt{d^2+4+8\sqrt{d-1}}}{2},\ d+2\right].
\end{equation}
As $\ell\to\infty$, both endpoints of $J_{d,\ell}$ tend to $d^2/(d-1)$.
Indeed, for fixed $x>4$, $p_{\ell-2}/p_{\ell-1}\to r(x)$ and
$1/p_{\ell-1}\to0$, where
\[
 r(x)=\frac{x-2-\sqrt{x(x-4)}}2.
\]
The limiting equation is $x=d(1+r(x))$, with unique solution $d^2/(d-1)$;
monotonicity brackets the finite roots around this solution. For $d=3$
the limit is $9/2$. The set of strict accumulation points is closed, so
$9/2\in\LimQ$ as well. For example,
\[
J_{3,2}=\left[\frac{5+\sqrt{13+8\sqrt2}}2,5\right],\qquad
J_{3,3}\approx[4.718680276,4.732050808].
\]
The approximations are for orientation; \eqref{eq:left}--\eqref{eq:right}
define the exact endpoints.

\subsection{Generalized Bethe decorations}
A word $\bb=(b_0,\ldots,b_{h-1})$ of positive integers specifies a finite
rooted generalized Bethe tree: every vertex at level $j$ has $b_j$ children,
and level $h$ consists of leaves. Let $\DD_{\bb}(G)$ identify the root of
one copy of this tree with \emph{every} vertex of $G$. The empty word means
no decoration. Compound graphs of this kind and their spectra have been
studied by Rojo and Medina \cite{RojoMedina}; the class itself is not new.

Define $F_{\bb}(\mu)=\lambda_{\max}(H_{\bb}(\mu))$, where $H_{\bb}$ is the
symmetric tridiagonal matrix with diagonal
\begin{equation}\label{eq:diag}
             (\mu+b_0,\ b_1+1,\ldots,b_{h-1}+1,\ 1)
\end{equation}
and off-diagonal entries $\sqrt{b_0},\ldots,\sqrt{b_{h-1}}$.
Set $F_{\varnothing}(\mu)=\mu$. For example,
\[
 H_{(q,b)}(\mu)=
 \begin{pmatrix}
 \mu+q&\sqrt q&0\\
 \sqrt q&b+1&\sqrt b\\
 0&\sqrt b&1
 \end{pmatrix}.
\]

\begin{proposition}[Exact decoration formula]\label{prop:bethe}
For every connected graph $G$,
\begin{equation}\label{eq:bethe}
       \rh(Q(\DD_{\bb}(G)))=F_{\bb}(\rh(Q(G))).
\end{equation}
The function $F_{\bb}$ is continuous and strictly increasing on
$[0,\infty)$. Hence for every $[u,v]\subseteq\LimQ$,
\[
                  [F_{\bb}(u),F_{\bb}(v)]\subseteq\LimQ.
\]
\end{proposition}
\begin{proof}
Let $z>0$ be a Perron vector of $Q(G)$ with eigenvalue $\mu$, and let
$y>0$ be a Perron vector of $H_{\bb}(\mu)$. The number of vertices at
level $j$ in one decoration is $N_0=1$ and
$N_j=b_0\cdots b_{j-1}$. At a level-$j$ vertex in the decoration rooted
at $v$, assign the value $z_vy_j/\sqrt{N_j}$. The eigenvalue equations
reduce to the equations for $H_{\bb}(\mu)y$, including the root equation
where $Q(G)z=\mu z$. This gives a positive eigenvector of the decorated
graph, proving~\eqref{eq:bethe} by Perron--Frobenius.

Continuity follows from the variational characterization. If
$\mu_2>\mu_1$, a unit top eigenvector $y$ at $\mu_1$ has $y_0>0$, and
\[
 F_{\bb}(\mu_2)\ge F_{\bb}(\mu_1)+(\mu_2-\mu_1)y_0^2
                      >F_{\bb}(\mu_1).
\]
Applying the exact formula to each strict approximating sequence proves
the interval assertion. Attaching trees preserves bipartiteness and
connectedness, so $L$ and $Q$ have the same radii throughout.
\end{proof}

\subsection{The finite covering argument}
\begin{proof}[Proof of Theorem~\ref{cor:LQ}]
For each of the 101 rows in Appendix~\ref{app:chain}, take the interval
\[
                    K_i=F_{\bb_i}(J_{d_i,\ell_i}).
\]
Theorem~\ref{thm:sub} and Proposition~\ref{prop:bethe} give $K_i\subseteq\LimQ$.
The accompanying exact certificate supplies rational inner intervals
$[u_i,v_i]\subseteq K_i$ satisfying
\begin{equation}\label{eq:coverchecks}
 u_1<29/5<v_1,\qquad u_i<v_{i-1}<v_i\ (2\le i\le101),
 \qquad v_{101}>6+2\sqrt5.
\end{equation}
Thus their union contains $[29/5,\Ga]$, and~\eqref{eq:oldray} completes
the proof. All assertions in~\eqref{eq:coverchecks} are checked by the
supplied verifier using integers and rational numbers only.
\end{proof}

The first interval is $F_{(1,1,1)}(J_{3,7})$, with an inner interval whose
endpoints are approximately
\[
             5.799661200535042,\qquad5.8004998873099805.
\]
The last is $F_{(1,6,4,3)}(J_{5,1})$; its certified right endpoint is
approximately $11.306678106894612$. The smallest consecutive overlap in
the chain is approximately $4.496717526157948\cdot10^{-7}$ and is
strictly positive as an exact rational number. The base degrees used are
$3,4,5$, the largest subdivision length is $8$, and the largest Bethe
word length is $9$. No translation or Cartesian-product rule is used.

The proof gives no covering of the whole interval $[9/2,29/5)$.
In particular, attaching Bethe trees to a graph does not remove its
cycles, and the input theorem is not restricted to trees. Results on
linear trees, such as those in \cite{BelardoOliveiraTrevisan}, therefore
remain a separate part of the Laplacian limit-point problem.

\section{Complementary tree constructions}\label{sec:trees}

\begin{lemma}[Strict growth for trees]\label{lem:tree-growth}
Suppose $\alpha\neq1/2$ and $T'$ is obtained from a tree $T$ by attaching at
least one new vertex.  Then
\[
 \rho_{B_\alpha}(T')>\rho_{B_\alpha}(T).
\]
\end{lemma}

\begin{proof}
By Lemma~\ref{lem:bipartite-switching}, it is enough to consider the
irreducible nonnegative matrix $H_{\alpha}=a_{\alpha}D+b_{\alpha}A$, where
$b_{\alpha}>0$.  Extend a positive Perron vector of $H_{\alpha}(T)$ by zeros
on the new vertices.  If $a_{\alpha}>0$, its Rayleigh quotient strictly
increases because the degree of an attachment vertex increases.  If
$a_{\alpha}=0$, the quotient is initially unchanged, but equality in the
spectral-radius comparison is impossible: the eigenvalue equation at a new
vertex would force a positive old coordinate to be zero.  Thus the
inequality is strict in both cases.
\end{proof}

At $\alpha=1/2$, one has $B_{1/2}=D/2$, and strict subtree monotonicity fails.

\subsection{Eigenvalue location on a path}

We use the standard diagonalization algorithm for a symmetric matrix whose
underlying graph is a tree \cite{JacobsTrevisan,HoppenJacobsTrevisan}.  Applied
to $B_{\alpha}(T)-\lambda I$, it initializes an edge by $c_{\alpha}$ and a
vertex $v$ by $a_{\alpha}d(v)-\lambda$.  Processing a child with current
diagonal value $z\neq0$ subtracts $c_{\alpha}^2/z$ from its parent.  By
Sylvester's law of inertia, the signs of the final diagonal entries count the
eigenvalues above, below, and equal to $\lambda$.

In the recurrence discussion assume $\alpha\ne1/2$, so that
$b_\alpha>0$.  The degree-matrix case is treated separately.
Along a pendant path, put $q_{\alpha}:=b_{\alpha}^2=c_{\alpha}^2$.  The
diagonal values satisfy
\begin{equation}\label{eq:path-recurrence}
 z_1=a_{\alpha}-\lambda,
 \qquad
 z_{j+1}=\phi_{\alpha,\lambda}(z_j),
 \qquad
 \phi_{\alpha,\lambda}(t)
 :=2a_{\alpha}-\lambda-\frac{q_{\alpha}}{t}.
\end{equation}
The discriminant of the fixed-point equation is
\begin{equation}\label{eq:discriminant}
 \Delta_{\alpha}(\lambda)
 :=(2a_{\alpha}-\lambda)^2-4q_{\alpha}.
\end{equation}
The relevant threshold is
\begin{equation}\label{eq:pi}
 \pi(\alpha):=2(a_{\alpha}+b_{\alpha})
 =\begin{cases}
   4-6\alpha,&0\leq\alpha\leq1/2,\\
   2\alpha,&1/2\leq\alpha\leq1.
  \end{cases}
\end{equation}
For $\lambda>\pi(\alpha)$, one has
$2a_{\alpha}-\lambda<-2b_{\alpha}$ and hence
$\Delta_{\alpha}(\lambda)>0$.  The two fixed points are
\begin{equation}\label{eq:fixed-points}
 \theta_{\alpha}(\lambda)
 =\frac{2a_{\alpha}-\lambda-\sqrt{\Delta_{\alpha}(\lambda)}}2,
 \qquad
 \theta'_{\alpha}(\lambda)
 =\frac{2a_{\alpha}-\lambda+\sqrt{\Delta_{\alpha}(\lambda)}}2.
\end{equation}
They satisfy
\[
 \theta_{\alpha}<-b_{\alpha}<\theta'_{\alpha}<0,
 \qquad
 \theta_{\alpha}\theta'_{\alpha}=q_{\alpha}.
\]
Moreover, $z_1<\theta_{\alpha}$ and the sequence in
\eqref{eq:path-recurrence} increases to the attracting fixed point
$\theta_{\alpha}$.  These are the standard hyperbolic properties of the
rational recurrence; see \cite{OliveiraTrevisanRecurrences}.
The convergence is uniform for $\lambda$ in any compact interval
$K\subset(\pi(\alpha),\infty)$.  Indeed, for $t\le\theta_\alpha(\lambda)$,
\[
0\le\phi'_{\alpha,\lambda}(t)
\le \max_{x\in K}\frac{q_\alpha}{\theta_\alpha(x)^2}<1.
\]
The initial values and the fixed points are bounded on $K$, so the
contraction estimate is uniform.

\subsection{The starlike trees \texorpdfstring{$T_{1,m,m}$}{T(1,m,m)}}

Let $T_{1,m,m}$ be the tree obtained by joining one path of length $1$ and two
paths of length $m$ at a common root.  The case $T_{1,3,3}$ supplies a uniform
lower bound which puts all subsequent calculations in the hyperbolic regime.

\begin{proposition}\label{prop:T133}
For every $\alpha\in[0,1]$,
\[
 \rho_{B_\alpha}(T_{1,3,3})\geq\pi(\alpha),
\]
and equality holds if and only if $\alpha=1$.
\end{proposition}

\begin{proof}
First suppose $b_{\alpha}>0$.  By bipartite switching, diagonalize
$H_{\alpha}(T_{1,3,3})-\pi(\alpha)I$, where
$H_{\alpha}=a_{\alpha}D+b_{\alpha}A$.  Write
$x=a_{\alpha}/b_{\alpha}$.  On either long arm the successive values are
\[
 z_1=-b_{\alpha}(x+2),\qquad
 z_2=-b_{\alpha}\frac{2x+3}{x+2},\qquad
 z_3=-b_{\alpha}\frac{3x+4}{2x+3},
\]
so all of them are negative.  After the short arm and the two long arms are
processed, the value at the root is
\begin{align*}
 \xi
 &=3a_{\alpha}-\pi(\alpha)
   -\frac{b_{\alpha}^2}{z_1}
   -2\frac{b_{\alpha}^2}{z_3}\\
 &=b_{\alpha}\left(
 x-2+\frac1{x+2}+2\frac{2x+3}{3x+4}\right)\\
 &=b_{\alpha}\frac{x(x+1)(3x+5)}{(x+2)(3x+4)}\geq0.
\end{align*}
It is positive when $a_{\alpha}>0$ and is zero when
$a_{\alpha}=0$, namely when $\alpha=1$.  The inertia theorem now gives the
claim.  If $b_{\alpha}=0$, then $\alpha=1/2$ and
$B_{1/2}=D/2$; hence
$\rho_{B_{1/2}}(T_{1,3,3})=3/2>1=\pi(1/2)$.
\end{proof}

\begin{remark}
There is no need for a closed expression for
$\rho_{B_\alpha}(T_{1,3,3})$.  Proposition~\ref{prop:T133} is exactly the
information required to guarantee that the path recurrence has two real
fixed points.  For $\alpha\ne1/2$ and $m\ge4$, strict tree growth gives
$\rho_{B_\alpha}(T_{1,m,m})>\pi(\alpha)$, including the adjacency case
where equality occurs at $m=3$.  Thus the limiting calculation takes
place strictly above the path threshold. 
\end{remark}

By Lemma~\ref{lem:tree-growth}, the sequence
$\rho_{B_\alpha}(T_{1,m,m})$ is strictly increasing in $m$ whenever
$\alpha\neq1/2$.  It is bounded, for example by the maximum row sum of the
switched nonnegative matrix, and therefore it converges.  At $\alpha=1/2$
the sequence is constantly $3/2$.

\begin{theorem}[The $T_{1,m,m}$ limit]\label{thm:tau0}
For each $\alpha\in[0,1]$, the limit
\[
 \tau_0(\alpha):=
 \lim_{m\to\infty}\rho_{B_\alpha}(T_{1,m,m})
\]
is the unique solution $\lambda>\pi(\alpha)$ of
\begin{equation}\label{eq:F0}
 F_{\alpha}(\lambda):=
 a_{\alpha}
 +\frac{q_{\alpha}}{\lambda-a_{\alpha}}
 -\sqrt{(\lambda-2a_{\alpha})^2-4q_{\alpha}}
 =0.
\end{equation}
If $\alpha\neq1/2$, then $\tau_0(\alpha)$ is a
$B_{\alpha}$-limit point.  At $\alpha=1/2$ the numerical limit is $3/2$ but
is not a limit point in the sense of Definition~\ref{def:limit-point}.
\end{theorem}

\begin{proof}
If $\alpha=1/2$, then $\rho_{B_{1/2}}(T_{1,m,m})=3/2$ and
$F_{1/2}(\lambda)=3/2-\lambda$ for $\lambda>1$, proving the assertion
directly.  Assume henceforth that $\alpha\ne1/2$ and $m\ge4$.
Let $\lambda_m=\rho_{B_\alpha}(T_{1,m,m})$.  Processing the three arms of
$T_{1,m,m}$ at $\lambda_m$ gives
\begin{equation}\label{eq:root-finite}
 3a_{\alpha}-\lambda_m
 -\frac{q_{\alpha}}{a_{\alpha}-\lambda_m}
 -2\frac{q_{\alpha}}{z_m(\lambda_m)}=0,
\end{equation}
where $z_m$ is generated by \eqref{eq:path-recurrence}.  The preceding
proposition and strict growth show that the limiting value is strictly above
$\pi(\alpha)$, including the exceptional equality at
$(\alpha,m)=(1,3)$.  The hyperbolic recurrence therefore converges uniformly
on a compact neighborhood of the limiting parameter, and
\[
 z_m(\lambda_m)\longrightarrow
 \theta_{\alpha}(\lambda),
 \qquad \lambda:=\lim_m\lambda_m.
\]
Taking limits in \eqref{eq:root-finite} and using
$q_{\alpha}/\theta_{\alpha}=\theta'_{\alpha}$ gives
\[
 3a_{\alpha}-\lambda
 -\frac{q_{\alpha}}{a_{\alpha}-\lambda}
 -2\theta'_{\alpha}(\lambda)=0.
\]
Since
$2\theta'_{\alpha}=2a_{\alpha}-\lambda
 +\sqrt{\Delta_{\alpha}(\lambda)}$,
this is exactly \eqref{eq:F0}.

For $\lambda>\pi(\alpha)$,
\[
 F'_{\alpha}(\lambda)
 =-\frac{q_{\alpha}}{(\lambda-a_{\alpha})^2}
  -\frac{\lambda-2a_{\alpha}}
  {\sqrt{(\lambda-2a_{\alpha})^2-4q_{\alpha}}}<0.
\]
Moreover,
$F_{\alpha}(\pi(\alpha)^+)>0$ and
$F_{\alpha}(\lambda)\to-\infty$ as $\lambda\to\infty$.
Thus \eqref{eq:F0} has exactly one solution in the required domain.  For
$\alpha\neq1/2$, the finite spectral radii are strictly increasing and hence
pairwise distinct; their limit is a limit point.  The assertion at
$\alpha=1/2$ follows from $B_{1/2}=D/2$.
\end{proof}

For symbolic computations one may square \eqref{eq:F0}, but the domain and
the sign in the unsquared equation must be retained.  With
$x=\lambda-a_{\alpha}$, the relevant root satisfies
\begin{equation}\label{eq:quartic}
 x^4-2a_{\alpha}x^3-4q_{\alpha}x^2
 -2a_{\alpha}q_{\alpha}x-q_{\alpha}^2=0,
 \qquad \lambda>\pi(\alpha),
\end{equation}
as well as \eqref{eq:F0}.  This restriction removes the spurious roots seen
in solving the squared equation without its domain restriction.

Some representative values are listed below.  We denote by
$\mathfrak h_L=4.382975767\ldots$ the classical $T_{1,m,m}$ Laplacian
limit.

\begin{table}[ht]
\centering
\caption{The path threshold, the $T_{1,m,m}$ limit, and the bound from weighted interpolation.}
\label{tab:thresholds}
\begin{tabular}{ccccc}
\toprule
$\alpha$ & $\pi(\alpha)$ & $\tau_0(\alpha)$ & $d_{\alpha}$ & $\Gamma(\alpha)$\\
\midrule
$0$ & $4$ & $\mathfrak h_L$ & $6$ & $6+2\sqrt5$\\
$1/4$ & $5/2$ & $2.825768143\ldots$ & $7$ & $21/4+\sqrt6$\\
$1/2$ & $1$ & $3/2$ & -- & --\\
$2/3$ & $4/3$ & $\mathfrak h_L/3$ & $6$ & $(6+2\sqrt5)/3$\\
$3/4$ & $3/2$ & $1.591081806\ldots$ & $5$ & $13/4$\\
$1$ & $2$ & $\sqrt{2+\sqrt5}$ & $4$ & $2\sqrt3$\\
\bottomrule
\end{tabular}
\end{table}

The column $\Gamma(\alpha)$ records only the bound from
Corollary~\ref{cor:halfline}. Theorem~\ref{cor:LQ} improves it to $29/5$
at $\alpha=0$, and to $29/15$ at $\alpha=2/3$ because $B_{2/3}=Q/3$.
The constructions of this section retain their separate interest as
tree realizations.

\subsection{Further starlike limit points}

\label{sec:starlike-gaps}

For positive integers $q_1\le\cdots\le q_m$, let
$S(q_1,\ldots,q_m)$ be the starlike tree whose arms have those lengths.
Assume throughout this subsection that $\alpha\neq1/2$, abbreviate
$a=a_{\alpha}$, $b=b_{\alpha}>0$, and $q=b^2$, and, for
$x>\pi(\alpha)$, define the branch responses
\begin{equation}\label{eq:branch-responses}
 R_1(x):=\frac{q}{x-a},\qquad
 R_{j+1}(x):=\frac{q}{x-2a-R_j(x)},
\end{equation}
and
\begin{equation}\label{eq:branch-response-infinity}
 R_{\infty}(x):=
 \frac{x-2a-\sqrt{(x-2a)^2-4q}}{2}.
\end{equation}
These are the positive Schur-complement contributions of, respectively, a
finite pendant path and a one-way infinite path.

\begin{lemma}\label{lem:branch-responses}
For every $x>\pi(\alpha)$,
\[
 0<R_1(x)<R_2(x)<\cdots<R_{\infty}(x)<b.
\]
Each $R_j$ and $R_{\infty}$ is strictly decreasing in $x$.  If
$\rho_{B_\alpha}(S(q_1,\ldots,q_m))>\pi(\alpha)$, then its spectral radius is
the unique solution in this range of
\begin{equation}\label{eq:starlike-root-equation}
 x-ma=\sum_{i=1}^m R_{q_i}(x).
\end{equation}
\end{lemma}

\begin{proof}
The map $u\mapsto q/(x-2a-u)$ is increasing on
$[0,R_{\infty}(x)]$ and fixes $R_{\infty}(x)$.  Moreover,
\[
 R_1(x)<R_{\infty}(x)
 \quad\Longleftrightarrow\quad
 x-a>\frac{q}{R_{\infty}(x)}
 =x-2a-R_{\infty}(x),
\]
which is immediate.  Iteration proves the first chain of inequalities.
Strict decrease in $x$ follows from \eqref{eq:branch-responses} by induction,
and directly from \eqref{eq:branch-response-infinity}.
For completeness, on a compact interval $K\subset(\pi(\alpha),\infty)$
the derivative of $u\mapsto q/(x-2a-u)$, for $0\le u\le R_\infty(x)$,
is at most
\[
\kappa_K:=\max_{x\in K}\frac{R_\infty(x)^2}{q}<1.
\]
Thus $0<R_\infty(x)-R_j(x)\le b\kappa_K^{j-1}$ uniformly on $K$.
This justifies the parameter-dependent limits used below.
Finally, eliminating
each arm from its leaf towards the root gives
\eqref{eq:starlike-root-equation}.  Its left side is strictly increasing and
its right side strictly decreasing, proving uniqueness.
\end{proof}

The next theorem gives both an absolute starlike band and a systematic source
of limit points within it.

\begin{theorem}[Starlike arms tending to infinity]
\label{thm:starlike-arms-infinity}
Let $m\geq5$ and $\alpha\neq1/2$.  Define
\begin{equation}\label{eq:tm-Lambda-m}
 t_m(\alpha):=
 \frac{\sqrt{(m-2)^2a_{\alpha}^2+4(m-1)b_{\alpha}^2}
       -(m-2)a_{\alpha}}
      {2(m-1)b_{\alpha}},
 \qquad
 \Lambda_m(\alpha):=ma_{\alpha}+mb_{\alpha}t_m(\alpha),
\end{equation}
and
\begin{equation}\label{eq:sigma-m}
 \sigma_m(\alpha):=
 \frac{(m+1)a_{\alpha}
 +\sqrt{(m-1)^2a_{\alpha}^2+4m b_{\alpha}^2}}{2}.
\end{equation}
Then:
\begin{enumerate}[label=\textup{(\roman*)},leftmargin=2.2em]
 \item $\Lambda_m(\alpha)$ is a $B_{\alpha}$-limit point, obtained from
       $S(n,\ldots,n)$ as $n\to\infty$.
 \item Every finite $m$-arm starlike tree satisfies
       \[
       \sigma_m(\alpha)
       \leq\rho_{B_\alpha}(S(q_1,\ldots,q_m))
       <\Lambda_m(\alpha),
       \]
       with equality on the left precisely for the star
       $S(1,\ldots,1)$.
 \item If $1\leq r<m$ and $q_1\leq\cdots\leq q_{m-r}$ are fixed, then
       \begin{equation}\label{eq:partial-starlike-limit}
       \ell_{\alpha}(q_1,\ldots,q_{m-r};r)
       :=\lim_{n\to\infty}
       \rho_{B_\alpha}\bigl(S(q_1,\ldots,q_{m-r},
                       \underbrace{n,\ldots,n}_{r})\bigr)
       \end{equation}
       exists, is a $B_{\alpha}$-limit point, and is the unique solution
       $x>\pi(\alpha)$ of
       \begin{equation}\label{eq:partial-starlike-equation}
       x-ma_{\alpha}
       =\sum_{i=1}^{m-r}R_{q_i}(x)+rR_{\infty}(x).
       \end{equation}
       It is strictly smaller than $\Lambda_m(\alpha)$.
\end{enumerate}
\end{theorem}

\begin{proof}
The quotient matrix of $aD+bA$ on the star $K_{1,m}$ has characteristic
factor
\[
 p_m(x)=(x-ma)(x-a)-mb^2,
\]
whose larger root is \eqref{eq:sigma-m}.  At
$x=\pi(\alpha)=2a+2b$, writing $y=a/b$ gives
\[
 \frac{p_m(\pi(\alpha))}{b^2}
 =-\bigl((m-2)y^2+2(m-3)y+(m-4)\bigr)<0.
\]
Thus $\sigma_m(\alpha)>\pi(\alpha)$ for $m\geq5$.

By Lemma~\ref{lem:tree-growth}, the spectral radii of
$S(n,\ldots,n)$ are strictly increasing and bounded.  Their root equations
are
\[
 x_n-ma=mR_n(x_n).
\]
Since the parameters stay in a compact subinterval above $\pi(\alpha)$,
$R_n$ converges uniformly there to $R_{\infty}$.  If $x_n\to x$ and
$t=R_{\infty}(x)/b$, the limiting root equation and the fixed-point equation
for $R_{\infty}$ give
\[
 x=ma+mbt=2a+b(t+t^{-1}).
\]
Consequently
\[
 (m-1)bt^2+(m-2)at-b=0.
\]
Its unique positive root is $t_m(\alpha)$ in
\eqref{eq:tm-Lambda-m}; it lies in $(0,1)$, so the limiting value is above
the path threshold.  This proves (i).

The star is a subtree of every $m$-arm starlike tree.  Conversely, every
finite $S(q_1,\ldots,q_m)$ is a proper subtree of $S(N,\ldots,N)$ for all
sufficiently large $N$.  Part (ii) now follows from strict tree growth and
part (i).

For (iii), increasing $n$ adds vertices to $r$ arms, so the finite spectral
radii are strictly increasing and bounded above by $\Lambda_m(\alpha)$.
Passing to the limit in their root equations gives
\eqref{eq:partial-starlike-equation}.  Uniqueness follows from
Lemma~\ref{lem:branch-responses}.  At $x=\Lambda_m(\alpha)$, the left side
minus the right side of \eqref{eq:partial-starlike-equation} equals
\[
 \sum_{i=1}^{m-r}
 \bigl(R_{\infty}(x)-R_{q_i}(x)\bigr)>0,
\]
so its root is strictly smaller than $\Lambda_m(\alpha)$.
\end{proof}

One particularly simple nested family is obtained by keeping $m-1$ equal
finite arms.

\begin{corollary}\label{cor:nested-starlike-points}
For $m\geq5$ and $\alpha\neq1/2$, put
\[
 \eta_{m,k}(\alpha)
 :=\ell_{\alpha}(\underbrace{k,\ldots,k}_{m-1};1),
 \qquad k\geq1.
\]
Then all $\eta_{m,k}(\alpha)$ are $B_{\alpha}$-limit points and
\[
 \sigma_m(\alpha)<\eta_{m,1}(\alpha)
 <\eta_{m,2}(\alpha)<\cdots<\Lambda_m(\alpha),
 \qquad
 \lim_{k\to\infty}\eta_{m,k}(\alpha)=\Lambda_m(\alpha).
\]
\end{corollary}

\begin{proof}
At a fixed $x>\pi(\alpha)$, Lemma~\ref{lem:branch-responses} gives
$R_{k+1}(x)>R_k(x)$.  Comparing the strictly increasing left-minus-right
functions in \eqref{eq:partial-starlike-equation} therefore shows that the
roots increase strictly with $k$.  The first one is strictly larger than the
star value by tree growth.  Finally, $R_k\to R_{\infty}$ uniformly on compact
subintervals above the path threshold, so the increasing limit of the roots
satisfies the all-infinite equation and equals $\Lambda_m(\alpha)$.
\end{proof}

For the Laplacian, $a_0=b_0=1$, and the formulas simplify to
\begin{equation}\label{eq:Laplacian-starlike-band}
 \sigma_m(0)=m+1,\qquad
 t_m(0)=\frac1{m-1},\qquad
 \Lambda_m(0)=m+\frac{m}{m-1}.
\end{equation}

For example, $\eta_{7,k}(0)\uparrow49/6$ and
$\eta_{9,k}(0)\uparrow81/8$, with
\[
8<\eta_{7,k}(0)<49/6<9,\qquad
10<\eta_{9,k}(0)<81/8<6+2\sqrt5.
\]
These are explicit tree realizations inside the half-line already covered
by Theorem~\ref{cor:LQ}. The next obstruction concerns only starlike
trees; it is not a gap in the full Laplacian limit-point set.

\begin{proposition}[An obstruction for starlike trees]
\label{prop:starlike-obstruction}
Neither $(49/6,9)$ nor $(81/8,6+2\sqrt5)$ contains a Laplacian spectral
radius of a finite starlike tree, or a limit of such radii.
\end{proposition}
\begin{proof}
A starlike with $m\ge5$ arms has index in
$[m+1,m+m/(m-1))$ by Theorem~\ref{thm:starlike-arms-infinity}.
For $m=5,6$ its index is at most $36/5<8$; for $m=7$ the upper limiting
bound is $49/6$; for $m=8$ the entire band is $[9,64/7]$; for $m=9$
the upper limiting bound is $81/8$; and for $m\ge10$ the lower bound
is at least $11$.  Trees with at most four arms have maximum degree at
most four and $\rho_L\le2\Delta\le8$.
These bounds persist under limits.  More explicitly, in a convergent
sequence with bounded indices the numbers of arms are bounded by the
lower bound $m+1$, so a subsequence has a fixed number of arms.
\end{proof}

\section{Relations with \texorpdfstring{$A_{\beta}$}{A-beta} and the deformed Laplacian}
\label{sec:relations}

\subsection{The deformed Laplacian}
For a real parameter $s$, let
\[
M_G(s)=I-sA(G)+s^2(D(G)-I).
\]
Thus $M_G(1)=L(G)$ and $M_G(-1)=Q(G)$.
Write $\mathcal L_M(s)$ for the strict limit-point set of its spectral
radii at fixed $s$.

\begin{proposition}[Explicit affine similarity on bipartite graphs]
\label{prop:deformed-equivalence}
Let $G$ be bipartite with a fixed bipartition $X\cup Y$, let $s\ne0$,
put $r=|s|$, and set $\zeta_r=(1+r)/(1+2r)$ and
$\beta_r=r/(1+r)$. Order the vertices by the bipartition and define
\[
 S_s=\begin{cases}
 I,&s<0,\\
 \operatorname{diag}(I_{|X|},-I_{|Y|}),&s>0.
 \end{cases}
\]
Then $S_s^{-1}=S_s$ and
\begin{equation}\label{eq:deformed-B-A-equivalence}
\begin{aligned}
 S_s M_G(s) S_s
 &=(1-r^2)I+r(1+2r)B_{\zeta_r}(G)\\
 &=(1-r^2)I+r(1+r)A_{\beta_r}(G).
\end{aligned}
\end{equation}
If $G$ has no isolated vertices, the spectral radius of $M_G(s)$ equals
its largest eigenvalue, and
\begin{equation}\label{eq:deformed-limit-map}
\rho(M_G(s))=1-r^2+r(1+2r)\rho_{B_{\zeta_r}}(G).
\end{equation}
Consequently the same positive affine map carries the limit-point sets
on any fixed class of connected bipartite graphs with at least two vertices
into one another.
\end{proposition}
\begin{proof}
Expand $M_G(s)=(1-r^2)I+r^2D-sA$. The definition of $S_s$ gives
$S_sDS_s=D$ and $S_s(-sA)S_s=rA$: for $s<0$ this is immediate,
whereas for $s>0$ it follows from bipartite switching. Therefore
\[
 S_sM_G(s)S_s=(1-r^2)I+r^2D+rA.
\]
The coefficient identities
\[
 r(1+2r)B_{\zeta_r}=r^2D+rA=r(1+r)A_{\beta_r}
\]
prove \eqref{eq:deformed-B-A-equivalence}. This is an orthogonal
similarity, expressed as an equality after conjugation. If there are no
isolated vertices, the conjugated matrix has diagonal
$1+r^2(d(v)-1)\ge1$ and nonnegative off-diagonal entries $rA_{uv}$.
Perron--Frobenius thus identifies its spectral radius with its largest
eigenvalue. Orthogonal similarity and the positive affine formula give
\eqref{eq:deformed-limit-map}, preserving convergence and pairwise
distinctness. The no-isolates qualification avoids relying on the
generally false assertion that an arbitrary scalar translation preserves
spectral radii.
\end{proof}

For $s\ge1$ one also has the direct identity, without bipartiteness,
\begin{equation}\label{eq:deformed-direct-positive}
M_G(s)=(1-s^2)I+s(2s-1)B_{(s-1)/(2s-1)}(G).
\end{equation}
For $s<0$, $S_s=I$, so the affine equality in
\eqref{eq:deformed-B-A-equivalence} holds directly for every graph,
without a bipartiteness assumption.  The bipartite construction removes the restriction
on the sign or size of $s$ in the interval corollary.

\begin{corollary}[Deformed-Laplacian intervals for every nonzero parameter]
\label{cor:deformed-intervals}
Let $s\ne0$, $r=|s|$, and $d\ge3$.  Then
\begin{equation}\label{eq:deformed-Id}
K_d(s)=\left[1-r^2+r^2d+2r\sqrt{d-1},\,
             1-r^2+(r^2+r)d\right]\subseteq\mathcal L_M(s).
\end{equation}
All points have connected bipartite realizations.  In particular, with
\[
d_s=\left\lceil(1+\sqrt{1+r})^2\right\rceil,\qquad
\Gamma_M(s)=1-r^2+r^2d_s+2r\sqrt{d_s-1},
\]
one has $[\Gamma_M(s),\infty)\subseteq\mathcal L_M(s)$.
At $s=0$, $M_G(0)=I$ and $\mathcal L_M(0)=\varnothing$.
\end{corollary}
\begin{proof}
Apply Theorem~\ref{thm:positive-interpolation} with $a=r^2$, $b=r$.
Its connected approximants eventually have at least two vertices, since
their indices are positive.  Proposition~\ref{prop:deformed-equivalence}
translates their indices by $1-r^2$, giving \eqref{eq:deformed-Id} and
preserving distinctness.  Adjacent intervals overlap exactly when
$r(d-2\sqrt d)\ge r^2$, which gives $d_s$ as in Section~\ref{sec:overlap}.
The assertion at zero is immediate.
\end{proof}

For $0<r<1$, $\beta_r<1/2$ and $\zeta_r>2/3$.  Thus the tree
$T_{1,m,m}$ and Shearer-caterpillar results in
\cite{Oliveira2027deformed} are affine images of the
corresponding $A_\beta$ constructions in \cite{OliveiraTrevisanAalpha}.
In particular, for every $s\ne0$,
\begin{equation}\label{eq:tau0-deformed-map}
\lim_m\rho(M_{T_{1,m,m}}(s))
=1-r^2+r(1+2r)\tau_0(\zeta_r).
\end{equation}
The other starlike limits transfer by the same formula.  These identities
separate the tree reparametrizations from the weighted interpolation theorem.

\subsection{The \texorpdfstring{$A_{\beta}$}{A-beta} scaling and Shearer caterpillars}

\begin{corollary}[$A_\beta$ intervals for the full nondegenerate range]
\label{cor:A-intervals}
For every $0\le\beta<1$ and integer $d\ge3$,
\[
[\beta d+2(1-\beta)\sqrt{d-1},d]
\subseteq\mathcal L_A^{\mathrm{bip}}(\beta).
\]
If $d_\beta=\left\lceil(1+1/\sqrt{1-\beta})^2\right\rceil$, then
\[
[\beta d_\beta+2(1-\beta)\sqrt{d_\beta-1},\infty)
\subseteq\mathcal L_A^{\mathrm{bip}}(\beta).
\]
For $\beta=1$ the matrix is $D$, whose strict limit-point set is empty.
\end{corollary}
\begin{proof}
Apply Theorem~\ref{thm:positive-interpolation} with $a=\beta$,
$b=1-\beta$, and use the overlap calculation with
$a/b=\beta/(1-\beta)$.
\end{proof}
The statement supplies a half-line for every $\beta<1$, including
$\beta\ge1/2$, beyond the parameter range of the caterpillar half-line
theorem in \cite{OliveiraTrevisanAalpha}.  Its sufficient lower endpoint
need not be optimal: at $\beta=0$ it is $2\sqrt3$, whereas Shearer's
bound is $\sqrt{2+\sqrt5}$.

Identity \eqref{eq:global-A-beta} has two consequences.  First, for every
$\alpha>1/2$,
\begin{equation}\label{eq:limit-set-scaling}
 \mathcal L_B(\alpha)
 =\alpha\,\mathcal L_A\!\left(\frac{1-\alpha}{\alpha}\right),
\end{equation}
where $\mathcal L_A(\beta)$ denotes the set of $A_{\beta}$-limit points.
Second, the direct $B_{\alpha}$ caterpillar recurrence is not a separate
dynamical system: after scaling, it is exactly the recurrence studied in
\cite{OliveiraTrevisanAalpha}.

For the caterpillar recurrence, assume $\alpha\ne1/2$, let
$\lambda>\pi(\alpha)$ and define
\begin{equation}\label{eq:drift}
 \delta_{\alpha}(\lambda)
 :=a_{\alpha}+\frac{q_{\alpha}}{\lambda-a_{\alpha}}>0.
\end{equation}
If $r_j$ leaves are attached to the $j$th spine vertex of a caterpillar,
the diagonalization outputs at the first and at each nonterminal spine vertex
obey
\begin{equation}\label{eq:caterpillar-correct}
 b_1=a_{\alpha}-\lambda+r_1\delta_{\alpha}(\lambda),
 \qquad
 b_{j+1}=\phi_{\alpha,\lambda}(b_j)
          +r_{j+1}\delta_{\alpha}(\lambda).
\end{equation}
At the terminal spine vertex the initial diagonal term is $a_{\alpha}$
rather than $2a_{\alpha}$; this standard end correction is exactly the one
appearing in the $A_{\beta}$ construction.
Thus one chooses $r_{j+1}$ so that
\[
 \theta'_{\alpha}-\delta_{\alpha}
 \leq b_{j+1}<\theta'_{\alpha}.
\]
When a quotient is an integer, the ordinary floor must be replaced by the
greatest integer strictly smaller than the quotient; equivalently,
$\lfloor x\rfloor_{<}=\lceil x\rceil-1$.  This convention enforces the strict upper inequality, including integer quotients.

The simplest convergence criterion requires
\begin{equation}\label{eq:simple-shearer-condition}
 \theta'_{\alpha}(\lambda)-\delta_{\alpha}(\lambda)>-b_{\alpha},
\end{equation}
because then $b_j^2<b_{\alpha}^2$ and the derivative products in the
eigenvalue-location criterion diverge.  As $\lambda\to\infty$,
\[
 \theta'_{\alpha}(\lambda)\longrightarrow0,
 \qquad
 \delta_{\alpha}(\lambda)\longrightarrow a_{\alpha}.
\]
Condition \eqref{eq:simple-shearer-condition} can therefore hold for all
sufficiently large $\lambda$ precisely when
$a_{\alpha}<b_{\alpha}$, that is, when $\alpha>2/3$.  This identifies exactly the parameter range of this sufficient criterion.

Let $\kappa_2(\alpha)$ be the unique solution, for $\alpha>2/3$, of
\begin{equation}\label{eq:kappa2}
 \delta_{\alpha}(\lambda)-\theta'_{\alpha}(\lambda)=b_{\alpha}.
\end{equation}
Equivalently,
\[
 \kappa_2(\alpha)
 =\alpha\tau_2\!\left(\frac{1-\alpha}{\alpha}\right),
\]
where $\tau_2$ is the threshold in
\cite[Theorem~4.4]{OliveiraTrevisanAalpha}.  Hence
\begin{equation}\label{eq:caterpillar-ray}
 [\kappa_2(\alpha),\infty)
 \subseteq\mathcal L_B(\alpha),
 \qquad \alpha>2/3,
\end{equation}
and every point in this half-line is obtained using caterpillar trees.
Near $\alpha=2/3$, the threshold in \eqref{eq:kappa2} diverges, whereas the
general threshold $\Gamma(\alpha)$ from Corollary~\ref{cor:halfline} remains
finite.  Thus the induced-subgraph theorem covers a regime not supplied by this caterpillar criterion.

The sharper low-parameter part of the $A_{\beta}$ result also transfers.
In particular, \cite[Corollary~5.3]{OliveiraTrevisanAalpha} gives the following.

\begin{corollary}[A sharp ray near adjacency]\label{cor:sharp-near-one}
If
\[
 \alpha>\frac{5+\sqrt5}{8}\approx 0.9045,
\]
then
\[
 [\tau_0(\alpha),\infty)\subseteq\mathcal L_B(\alpha).
\]
\end{corollary}

\begin{proof}
Under $\beta=(1-\alpha)/\alpha$, the displayed condition is equivalent to
$0\leq\beta<1-2\sqrt5/5$.  Oliveira and Trevisan \cite{OliveiraTrevisanAalpha} prove that every number at
least their $T_{1,m,m}$ limit is an $A_{\beta}$-limit point throughout this
range.  Multiply their limit points by $\alpha$ and use
\eqref{eq:limit-set-scaling}.  The scaled starlike limit is precisely
$\tau_0(\alpha)$ by Theorem~\ref{thm:tau0}.  
\end{proof}

\subsection{The affine family generated by \texorpdfstring{$I$, $D$, and $A$}{I, D, and A}}
\label{subsec:general-pencil}
The families above are instances of
$\Psi_{u,v,w}(G)=uI+vD(G)+wA(G)$.  Their equivalences depend both on the
coefficient domain and on the graph class.

\begin{proposition}\label{prop:general-pencil}
If $v\ge0$, $w\ge-v$, and $(v,w)\ne(0,0)$, then
\begin{equation}\label{eq:general-pencil-B}
\Psi_{u,v,w}=uI+(2v+w)B_{(v+w)/(2v+w)},
\end{equation}
where the scale is positive and the $B$-parameter lies in $[0,1]$.
These conditions describe exactly the coefficient pairs obtained from
positive multiples of $B_\alpha$, $0\le\alpha\le1$.
If $G$ is bipartite and $v\ge0$, switching additionally gives
\[
\Psi_{u,v,w}(G)\sim uI+(2v+|w|)
B_{(v+|w|)/(2v+|w|)}(G)
\]
whenever $(v,w)\ne(0,0)$.
\end{proposition}
\begin{proof}
Writing $tB_\alpha=t(1-\alpha)D+t(2\alpha-1)A$, coefficient comparison
gives $t=2v+w$ and $\alpha=(v+w)/(2v+w)$.  The stated inequalities
are precisely $t>0$ and $0\le\alpha\le1$.  On bipartite graphs the
coefficient of $A$ can be changed from $w$ to $|w|$ by switching.
\end{proof}

Without the coefficient or bipartiteness restrictions, the equivalence
need not hold uniformly over graphs.  For example $D-2A$ cannot equal
$uI+tB_\alpha$ with fixed coefficients for every graph: comparison of
$D$ and $A$ would force $t=0$.  On regular graphs the relation $D=dI$
reduces the dimension of the family, but a negative adjacency coefficient
selects the least adjacency eigenvalue for the largest eigenvalue.
Specifically, for a $d$-regular graph and $w<0$,
\[
\lambda_{\max}(\Psi_{u,v,w})=u+vd+w\lambda_n(A),
\]
whereas for any real $w$,
\[
\rho(\Psi_{u,v,w})=
\max\{|u+vd+wd|,\ |u+vd+w\lambda_n(A)|\}.
\]
Thus a scalar translation preserves the affine law for ordered
eigenvalues, but not in general for spectral radii.  For example
$-I-A(K_3)$ has largest eigenvalue $0$ and spectral radius $3$.

\section{Concluding remarks and open directions}
The weighted interpolation theorem supplies explicit intervals and a
half-line for every $B_\alpha$ except the discrete case $\alpha=1/2$,
with connected bipartite realizations throughout. The same argument
gives half-lines for every $A_\beta$ with $\beta<1$ and every nonzero
deformed-Laplacian parameter. Subdivision and generalized Bethe
decoration improve the common Laplacian and signless Laplacian
half-line to $[29/5,\infty)$, with a finite exact covering certificate.

The tree formulas complement this result by describing explicit
sequences and the limitations of the starlike subclass. The former
intervals $(49/6,9)$ and $(81/8,6+2\sqrt5)$ are covered by
Theorem~\ref{cor:LQ}; they remain inaccessible to starlike trees alone.
The following questions concern what is not settled by these arguments.
\begin{enumerate}[label=\textup{(P\arabic*)},leftmargin=2.2em]
\item Determine the least onset of a full half-line in
$\mathcal L_B(\alpha)$, and compare it with the corresponding onset for
connected bipartite graphs. In the Laplacian case, can the bound $29/5$
be lowered, and can the whole interval $[9/2,29/5)$ be covered?
\item Determine which additional intervals admit realizations by trees,
or specifically by linear trees. The bipartite covering theorem does
not settle this restriction.
\item Determine the first $B_\alpha$-limit points below $\tau_0(\alpha)$
and the parameters for which the full half-line starting at
$\tau_0(\alpha)$ is filled.
\item Quantify the cost of constructing approximating graphs and improve,
where possible, the $O(N^{-1/2})$ accuracy of Remark~\ref{rem:rate}.
\end{enumerate}

\section*{Acknowledgments}
  Francesco Belardo and Callum Huntington acknowledge partial support
  from the Gruppo Nazionale per le Strutture Algebriche, Geometriche e le loro
  Applicazioni (GNSAGA) of the Istituto Nazionale di Alta Matematica
  (INdAM), Italy.
  The work of Elismar R. Oliveira and Vilmar Trevisan was partially funded
  by CNPq (Brazil), grant 408180/2023-4.
  Vilmar Trevisan also received support from CNPq (Brazil),
  grant 308774/2025-6.

\section*{Declaration of generative AI use}
During the preparation of this manuscript, the authors used ChatGPT (OpenAI) to assist
with specific task, e.g., generating an initial outline and brainstorming arguments. The
authors reviewed, edited, and verified all generated content, and take full responsibility
for the accuracy and originality of the final submission.

\clearpage
\appendix
\section{The 101-interval certificate}\label{app:chain}
Each row below specifies $K_i=F_{\bb_i}(J_{d_i,\ell_i})$. The companion
JSON file records exact rational
brackets and inner intervals $[u_i,v_i]\subseteq K_i$. The last two
columns show $u_i,v_i$ rounded to nine decimal places; those displayed
decimals are not used in the proof. The empty word $\varnothing$ means
no decoration.
\subsection*{Exact verification}
The JSON file records rational brackets for each endpoint at two stages.
\begin{enumerate}
\item For the base interval, put
$P(x)=(x-d)p_{\ell-1}(x)-dp_{\ell-2}(x)$. To test the sign of
$P(x)-2\sqrt{d-1}$, first test whether $P(x)<0$; otherwise test the sign
of $P(x)^2-4(d-1)$. The upper endpoint uses $P(x)-d$. Uniqueness proved
in Theorem~\ref{thm:sub} makes these exact root brackets.
\item Use the \emph{upper} bound of the lower base endpoint and the
\emph{lower} bound of the upper base endpoint. They form an inner base
interval. Monotonicity of $F_{\bb}$ preserves this inclusion.
\item For rational $\mu$, locate $F_{\bb}(\mu)$ by testing whether
$xI-H_{\bb}(\mu)$ is positive definite. Its successive $LDL^T$ pivots are
\[
 t_0=x-\mu-b_0,\qquad
 t_j=x-c_j-b_{j-1}/t_{j-1}\quad(1\le j\le h),
\]
where $c_j=b_j+1$ for $j<h$ and $c_h=1$. All pivots must be strictly
positive; a nonpositive pivot ends the test. Since only squared
off-diagonal entries occur, every operation is rational.
\item Round inward again at the two image endpoints, check all strict
overlaps, and check $v_{101}>6$ and $(v_{101}-6)^2>20$. The latter two
inequalities establish $v_{101}>\Ga$ without evaluating $\sqrt5$.
\end{enumerate}
The generator used 64 bisections for each bracket. The independent
verifier checks the recorded brackets themselves; it does not trust the
generator, floating-point eigenvalues, rounded table entries or a plotted
appearance of overlap. Python's standard library is sufficient:
\begin{verbatim}
python3 verify_bethe_cover.py Bethe_interval_certificate.json
\end{verbatim}
The program checks every row, all parameter restrictions, all root and
eigenvalue brackets, all inclusions and overlaps, the claimed minimum
overlap and the final endpoint. Its checks remain active under
\texttt{python3 -O}.

\clearpage
\begingroup\small
\setlength{\tabcolsep}{4pt}
\renewcommand{\arraystretch}{1.00}
\begin{longtable}{rrrlrr}
\toprule
No. & $d$ & $\ell$ & $\mathbf b$ & Inner left & Inner right\\\midrule
\endfirsthead
\toprule
No. & $d$ & $\ell$ & $\mathbf b$ & Inner left & Inner right\\\midrule
\endhead
\midrule\multicolumn{6}{r}{Continued on next page}\\\endfoot
\bottomrule\endlastfoot
1 & 3 & 7 & $(1,1,1)$ & 5.799661201 & 5.800499887\\
2 & 3 & 7 & $(1,1,1,1,2,3)$ & 5.800486144 & 5.801322277\\
3 & 3 & 7 & $(1,1,1,2,1,2,3)$ & 5.801318372 & 5.802153457\\
4 & 3 & 7 & $(1,1,1,2,2,2,1,3,1)$ & 5.802082971 & 5.802915367\\
5 & 3 & 7 & $(1,1,1,1,3,2,1,2,2)$ & 5.802673577 & 5.803481771\\
6 & 3 & 8 & $(1,1,2)$ & 5.803468633 & 5.803896879\\
7 & 3 & 7 & $(1,1,1,1,3,2,1,3)$ & 5.803720404 & 5.804500546\\
8 & 3 & 7 & $(1,1,1,3)$ & 5.804188761 & 5.805017309\\
9 & 3 & 7 & $(1,1,1,1,3,2,1,3,1)$ & 5.804972366 & 5.805707967\\
10 & 3 & 7 & $(1,1,1,2,2,3)$ & 5.805520334 & 5.806313968\\
11 & 3 & 7 & $(1,1,1,3,1,1,3,1,2)$ & 5.806303230 & 5.807122416\\
12 & 3 & 7 & $(1,1,1,3,1,2,2,2,2)$ & 5.807078825 & 5.807892713\\
13 & 3 & 7 & $(1,1,1,2,2,2,3)$ & 5.807546041 & 5.808255221\\
14 & 3 & 6 & $(1,1)$ & 5.808118379 & 5.809726712\\
15 & 3 & 7 & $(1,1,1,3,1,3,1)$ & 5.809719389 & 5.810511208\\
16 & 3 & 7 & $(1,1,2)$ & 5.810408947 & 5.811235363\\
17 & 3 & 7 & $(1,1,1,2,3,1,1,3,1)$ & 5.810944125 & 5.811683010\\
18 & 3 & 8 & $(1,1,2,1,2,3,1,2)$ & 5.811490305 & 5.811879726\\
19 & 3 & 8 & $(1,1,2,2,1,1,2,2,2)$ & 5.811875855 & 5.812294517\\
20 & 3 & 8 & $(1,1,2,2,1,2,2,2,1)$ & 5.812286325 & 5.812703985\\
21 & 3 & 7 & $(1,1,2,1,1,2,2,2,2)$ & 5.812702386 & 5.813524496\\
22 & 3 & 6 & $(1,1,1)$ & 5.813171732 & 5.814771319\\
23 & 3 & 6 & $(1,1,1,2,1,3)$ & 5.814759035 & 5.816352904\\
24 & 3 & 6 & $(1,1,1,2,2,2,2,2,2)$ & 5.815862389 & 5.817444929\\
25 & 3 & 6 & $(1,1,1,1,3,2,1,3,1)$ & 5.817175711 & 5.818659668\\
26 & 3 & 6 & $(1,1,1,2,2,3)$ & 5.818376845 & 5.819908046\\
27 & 3 & 4 & $(1)$ & 5.819150394 & 5.825276019\\
28 & 3 & 6 & $(1,1,1,2,3,1,2,2)$ & 5.825219295 & 5.826599100\\
29 & 3 & 6 & $(1,1,2,1,1,3,2,1)$ & 5.826579750 & 5.828137027\\
30 & 3 & 6 & $(1,1,2,1,2,2,2,3)$ & 5.828079218 & 5.829591573\\
31 & 3 & 1 & $\varnothing$ & 5.828427125 & 6.000000000\\
32 & 3 & 3 & $(1,1,2,1,3,2,1,3)$ & 5.999969580 & 6.012004202\\
33 & 3 & 3 & $(1,1,2,2,2,3,1,2,2)$ & 6.012001347 & 6.022873182\\
34 & 3 & 3 & $(1,1,2,3,1,1,3,2,1)$ & 6.022749407 & 6.033651193\\
35 & 3 & 3 & $(1,1,3,1,2,2,3)$ & 6.033177882 & 6.044327148\\
36 & 3 & 3 & $(1,1,3,1,3,1,3,2,2)$ & 6.043962922 & 6.053988946\\
37 & 3 & 3 & $(1,1,3,1,3,2,2,1,1)$ & 6.053846575 & 6.063042273\\
38 & 3 & 3 & $(1,1,3,2,1,2,2,2,3)$ & 6.062900497 & 6.072706984\\
39 & 3 & 3 & $(1,1,3,2,1,3,2,1,3)$ & 6.070826729 & 6.079925776\\
40 & 3 & 3 & $(1,1,3,2,2,1,3)$ & 6.079857133 & 6.088649590\\
41 & 3 & 3 & $(1,1,3,2,2,2,1,1,4)$ & 6.088267739 & 6.096018964\\
42 & 3 & 3 & $(1,1,3,2,2,2,2,2,3)$ & 6.093622842 & 6.101284939\\
43 & 3 & 3 & $(1,1,3,2,2,2,2,3,1)$ & 6.099892452 & 6.106890060\\
44 & 3 & 3 & $(1,2,1,1,3,2,1,3,2)$ & 6.106803197 & 6.118151688\\
45 & 3 & 3 & $(1,2,1,3,1,3,2)$ & 6.118044877 & 6.128957106\\
46 & 3 & 3 & $(1,2,1,3,2,2,2,1,3)$ & 6.128831025 & 6.138999377\\
47 & 3 & 3 & $(1,2,1,3,2,2,3)$ & 6.137603488 & 6.146691328\\
48 & 3 & 3 & $(1,2,2,1,3,2,2,2,1)$ & 6.146485469 & 6.156875799\\
49 & 3 & 3 & $(1,2,2,2,2,2,2)$ & 6.156806724 & 6.167058333\\
50 & 3 & 2 & $(1)$ & 6.159272418 & 6.192582404\\
51 & 3 & 4 & $(1,2,3,1,2,3,1,3,1)$ & 6.192498948 & 6.195943759\\
52 & 3 & 3 & $(1,2,1,4)$ & 6.193985457 & 6.200168466\\
53 & 3 & 3 & $(1,2,2,3,1,2,2,1,3)$ & 6.200167988 & 6.208779890\\
54 & 3 & 3 & $(1,2,2,3,1,3,1,3,2)$ & 6.208559523 & 6.216621244\\
55 & 3 & 2 & $(1,1)$ & 6.214062270 & 6.246551365\\
56 & 3 & 2 & $(1,1,3,1,2,2,3,1,3)$ & 6.246463150 & 6.277185095\\
57 & 3 & 2 & $(1,1,3,2,2,3,1,2,3)$ & 6.276826874 & 6.303020612\\
58 & 3 & 2 & $(1,1,3,2,3,1,3,1,3)$ & 6.300809336 & 6.322641473\\
59 & 3 & 2 & $(1,2,1,2,3,1,3)$ & 6.322620779 & 6.352629019\\
60 & 3 & 2 & $(1,2,2,2,2,1,2,2,3)$ & 6.352628569 & 6.381273503\\
61 & 3 & 2 & $(1,2,2,3,1,3,2,3)$ & 6.380743294 & 6.406234581\\
62 & 3 & 2 & $(1,2,2,3,2,2,3)$ & 6.404722294 & 6.426892821\\
63 & 3 & 2 & $(1,2,3,1,3,2,2,3,2)$ & 6.426350543 & 6.450125092\\
64 & 3 & 2 & $(1,2,3,2,1)$ & 6.449860310 & 6.472531761\\
65 & 3 & 2 & $(1,2,3,2,2,2,3,2,3)$ & 6.472142176 & 6.491895329\\
66 & 3 & 2 & $(1,3)$ & 6.482238131 & 6.509040333\\
67 & 3 & 2 & $(1,2,3,2,3,1,3,3)$ & 6.508166307 & 6.524261692\\
68 & 3 & 2 & $(1,3,1,2,3,1,4)$ & 6.524148421 & 6.549094148\\
69 & 3 & 2 & $(1,3,1,3,3)$ & 6.548485431 & 6.571106167\\
70 & 3 & 2 & $(1,3,1,3,3,1,4)$ & 6.568373455 & 6.587156946\\
71 & 3 & 2 & $(1,3,2,1,3,2,3,2)$ & 6.587145395 & 6.609688509\\
72 & 3 & 2 & $(1,3,2,2,2,3,2,2,4)$ & 6.609512119 & 6.629587805\\
73 & 3 & 2 & $(1,3,2,3)$ & 6.628899387 & 6.649255919\\
74 & 3 & 2 & $(1,3,2,3,1,3,3,2,2)$ & 6.649212608 & 6.667088540\\
75 & 3 & 2 & $(1,3,2,3,2,2,3,2,3)$ & 6.666221560 & 6.683281550\\
76 & 3 & 2 & $(1,3,2,3,2,3,2,2,3)$ & 6.682908978 & 6.698201277\\
77 & 4 & 4 & $(1,1,4,1,2,3,2,2,1)$ & 6.698196798 & 6.706340288\\
78 & 3 & 2 & $(1,3,2,3,2,3,2,3,3)$ & 6.704411073 & 6.715963512\\
79 & 4 & 3 & $(1)$ & 6.707872203 & 6.735893565\\
80 & 4 & 4 & $(1,2,2,3)$ & 6.734193101 & 6.743200947\\
81 & 3 & 2 & $(1,3,3,1)$ & 6.737502289 & 6.754425712\\
82 & 4 & 3 & $(1,1,1,3)$ & 6.753836379 & 6.781292490\\
83 & 4 & 3 & $(1,1,3,3)$ & 6.777347068 & 6.803505485\\
84 & 4 & 3 & $(1,1,4,1)$ & 6.800587161 & 6.825400140\\
85 & 4 & 3 & $(1,2)$ & 6.820306162 & 6.846874841\\
86 & 4 & 3 & $(1,2,2,2)$ & 6.842203728 & 6.868091652\\
87 & 5 & 2 & $\varnothing$ & 6.854101966 & 7.000000000\\
88 & 3 & 1 & $(1)$ & 6.995226494 & 7.162277660\\
89 & 3 & 1 & $(1,2,3,3)$ & 7.145512149 & 7.294347122\\
90 & 3 & 1 & $(1,3,3,1)$ & 7.287641072 & 7.421919844\\
91 & 3 & 1 & $(1,4)$ & 7.393930223 & 7.523476271\\
92 & 4 & 1 & $\varnothing$ & 7.464101615 & 8.000000000\\
93 & 5 & 2 & $(1)$ & 7.997019952 & 8.140054945\\
94 & 3 & 1 & $(2)$ & 8.109731692 & 8.274917218\\
95 & 3 & 1 & $(2,2,4,1)$ & 8.272018779 & 8.424834874\\
96 & 3 & 1 & $(2,3,4)$ & 8.420303241 & 8.560159706\\
97 & 3 & 1 & $(2,4,2,4)$ & 8.558009255 & 8.689853903\\
98 & 4 & 1 & $(1,2,5,2)$ & 8.683935365 & 9.187442196\\
99 & 5 & 1 & $\varnothing$ & 9.000000000 & 10.000000000\\
100 & 4 & 1 & $(2,4,5,2)$ & 9.995049618 & 10.455282688\\
101 & 5 & 1 & $(1,6,4,3)$ & 10.453686744 & 11.306678107\\
\end{longtable}
\endgroup


{\small
\begin{thebibliography}{99}

\bibitem{AlonWei}
N.~Alon and F.~Wei,
\newblock The limit points of the top and bottom eigenvalues of regular graphs,
\newblock arXiv:2304.01281v2, 2023.

\bibitem{BelardoOliveiraTrevisan}
F.~Belardo, E.~R.~Oliveira, and V.~Trevisan,
\newblock Limit points of (signless) Laplacian spectral radii of linear trees,
\newblock \emph{Appl. Math. Comput.} 477 (2024), 128819.
\newblock \url{https://doi.org/10.1016/j.amc.2024.128819}.


\bibitem{Diaz2026spectral}
R.~C.~D\'{i}az, E.~R.~Oliveira, and V.~Trevisan.
\newblock Spectral properties of the deformed Laplacian matrix of trees and {$H$}-join graphs.
\newblock \emph{SIAM J. Matrix Anal. Appl.} 47 (2026), no.~3, 1530--1550.
\newblock \href{https://doi.org/10.1137/25M1759823}{doi:10.1137/25M1759823}.


\bibitem{Friedman}
J.~Friedman,
\newblock A proof of Alon's second eigenvalue conjecture and related problems,
\newblock \emph{Mem. Amer. Math. Soc.} 195 (2008), no.~910.

\bibitem{Guo2008}
J.-M.~Guo,
\newblock On limit points of Laplacian spectral radii of graphs,
\newblock \emph{Linear Algebra Appl.} 429 (2008), 1705--1718.

\bibitem{Hoffman1972}
A.~J.~Hoffman,
\newblock On limit points of spectral radii of non-negative symmetric integral matrices,
\newblock in \emph{Graph Theory and Applications}, Lecture Notes in Math. 303,
Springer, 1972, pp.~165--172.

\bibitem{HoppenJacobsTrevisan}
C.~Hoppen, D.~P.~Jacobs, and V.~Trevisan,
\newblock \emph{Locating Eigenvalues in Graphs: Algorithms and Applications},
\newblock Springer, 2022.

\bibitem{JacobsTrevisan}
D.~P.~Jacobs and V.~Trevisan,
\newblock Locating the eigenvalues of trees,
\newblock \emph{Linear Algebra Appl.} 434 (2011), 81--88.

\bibitem{MSS}
A.~W.~Marcus, D.~A.~Spielman, and N.~Srivastava,
\newblock Interlacing families I: Bipartite Ramanujan graphs of all degrees,
\newblock \emph{Ann. of Math.} 182 (2015), 307--325.
\newblock \url{https://arxiv.org/abs/1304.4132}.


\bibitem{Nikiforov}
V.~Nikiforov,
\newblock Merging the $A$- and $Q$-spectral theories,
\newblock \emph{Appl. Anal. Discrete Math.} 11 (2017), 81--107.

\bibitem{OliveiraTrevisanRecurrences}
E.~R.~Oliveira and V.~Trevisan,
\newblock Applications of rational difference equations to spectral graph theory,
\newblock \emph{J. Difference Equ. Appl.} 27 (2021), 1024--1051.

\bibitem{OliveiraTrevisanAalpha}
E.~R.~Oliveira and V.~Trevisan,
\newblock Limit points of $A_{\alpha}$-matrices of graphs,
\newblock \emph{Linear Algebra Appl.} 728 (2026), 1--25.
\newblock \url{https://doi.org/10.1016/j.laa.2025.08.014}.


\bibitem{Oliveira2027deformed}
	E.~R.~Oliveira, J.~Szutkoski, and V.~Trevisan.
	\newblock Distribution of deformed Laplacian limit points.
	\newblock \emph{Linear Algebra Appl.} 753 (2027), 169--197.
	\newblock \href{https://doi.org/10.1016/j.laa.2026.08.032}{doi:10.1016/j.laa.2026.08.032}.


\bibitem{RojoMedina}
O.~Rojo and L.~Medina,
\newblock Spectra of weighted compound graphs of generalized Bethe trees,
\newblock \emph{Electron. J. Linear Algebra} 18 (2009), 30--57.
\newblock \url{https://doi.org/10.13001/1081-3810.1292}.

\bibitem{SamantaDeepshikhaDas}
A.~Samanta, Deepshikha, and K.~C.~Das,
\newblock Unifying adjacency, Laplacian, and signless Laplacian theories,
\newblock \emph{Ars Math. Contemp.} 24 (2024), Paper P4.10.
\newblock \url{https://doi.org/10.26493/1855-3974.3163.6hw}.

\bibitem{Shearer1989}
J.~B.~Shearer,
\newblock On the distribution of the maximum eigenvalue of graphs,
\newblock \emph{Linear Algebra Appl.} 114/115 (1989), 17--20.

\end{thebibliography}


}
\end{document}